\documentclass[graybox]{svmult}

\usepackage{mathptmx}       
\usepackage{helvet}         
\usepackage{courier}        
\usepackage{type1cm}        
\usepackage{makeidx}         
\usepackage{graphicx}        
\usepackage{multicol}        
\usepackage[bottom]{footmisc}
\usepackage{pstricks}

\usepackage{amsmath,amssymb, graphicx, nicefrac}
\usepackage{psfrag,url}

\usepackage{algpseudocode,algorithm} 
\floatname{algorithm}{Pseudocode} 

\usepackage{color}
\definecolor{darkred}{RGB}{139,0,0}
\definecolor{darkgreen}{RGB}{0,100,0}
\definecolor{darkmagenta}{RGB}{139,0,139}
\definecolor{darkpurple}{RGB}{110,0,180}
\definecolor{darkblue}{RGB}{40,0,200}
\definecolor{darkorange}{RGB}{255,140,0}

\newcommand{\beqn}{\begin{equation}}
\newcommand{\eeqn}{\end{equation}}

\newcommand{\fa}{{\mathfrak a}}

\def\R{\mathbb{R}}
\def\N{\mathbb{N}}

\newcommand{\cX}{\mathcal{X}}
\newcommand{\bcX}{\mathcal{X}}
\newcommand{\bcY}{\mathcal{Y}}
\newcommand{\cY}{\mathcal{Y}}

\newcommand{\bsbeta}{\boldsymbol{\beta}}
\newcommand{\bsgamma}{{\boldsymbol{\gamma}}}

\newcommand{\bsnu}{{\boldsymbol{\nu}}}

\newcommand{\bsy}{{\boldsymbol{y}}}

\newcommand{\bsu}{{\boldsymbol{u}}}

\newcommand{\bbR}{\mathbb{R}}

\newcommand{\bbN}{\mathbb{N}}

\newcommand{\calO}{\mathcal{O}}

\newcommand{\cL}{\mathcal{L}}

\newcommand{\mask}[1]{{}}

\numberwithin{equation}{section}

  \newtheorem{assumption}{Assumption}

\newcommand{\be}{\begin{equation}}
\newcommand{\ee}{\end{equation}}
\newcommand{\ba}{\begin{array}}
\newcommand{\ea}{\end{array}}
\newcommand{\beas}{\begin{eqnarray*}}
\newcommand{\eeas}{\end{eqnarray*}}
\newcommand{\bea}{\begin{eqnarray}}
\newcommand{\eea}{\end{eqnarray}}

\begin{document}

\title*{Higher Order Quasi Monte-Carlo Integration in Uncertainty Quantification}
\titlerunning{Higher Order QMC for UQ}

\author{Josef~Dick, Quoc~T.~Le Gia, Christoph~Schwab}

\institute{
Josef Dick and Quoc T. Le Gia \at 
School of Mathematics and Statistics, UNSW Australia, Sydney, Australia.
\email{josef.dick@unsw.edu.au, qlegia@unsw.edu.au}
\and
Christoph Schwab \at Seminar for Applied Mathematics, ETH, 8092 Z\"urich, Switzerland, 
\email{schwab@math.ethz.ch}
}
%
%
\maketitle

\abstract{
We review recent results on dimension-robust higher order convergence rates of Quasi-Monte Carlo Petrov-Galerkin approximations for response functionals of infinite-dimensional, parametric operator equations which arise in computational uncertainty quantification. 
}
\section{Introduction}
\label{sec:intro}
Computational uncertainty quantification (UQ) for
partial differential equations (PDEs)
with uncertain distributed input data gives rise, 
upon uncertainty parametrization, 
to the task of numerical solution of 
parametric, deterministic operator equations. 
Due to the distributed nature of uncertain inputs,
the number of parameters (and, hence, the
dimension of the parameter spaces) 
in such UQ problems is infinite.
The computation of response statistics corresponding to 
distributed uncertain inputs of PDEs involves, in addition,
\emph{numerical quadrature} of all possible `uncertain
scenarios', i.e., over the entire, infinite-dimensional
parameter space.

This has lead to the widespread use of sampling,
in particular Monte-Carlo (MC) and Markov-Chain Monte-Carlo (MCMC)
methods, in the numerical treatment of these problems:
MC methods afford convergence \emph{rates} which are
independent of the parameter dimension if the variance of the integrand
can be bounded independently of the dimension (the computational
\emph{work} of MC methods, of course, increases linearly with 
the space dimension).
This \emph{dimension robustness} of MC methods is purchased
at the cost of low order: the convergence rate of simple MC methods
is, generically, limited to $1/2$: variance reduction and
other devices can only reduced the constant, not the rate
in the convergence bounds.
At the same time, however, the \emph{parametric regularity}
required of integrand functions by MC methods is very moderate: 
mere square integrability with respect to a probability 
measure on the parameter space of the integrand functions
is needed, and point evaluations of the integrand functions
must be defined.
In UQ for problems whose solutions exhibit 
propagation of singularities
(as, eg., nonlinear hyperbolic conservation
laws with random inputs, see eg. 
\cite{MScSu2013,MScSu2014} and the references there), 
this kind of regularity is the best that can 
generally be expected. 
In other applications, the parametric dependence
of the response maps is considerably more regular:
the solutions' dependence on the parameters is, in fact, 
\emph{analytic}.
This observation has been the basis for
the widespread use of spectral- and polynomial chaos based
numerical methods for approximating the parameter dependence
in such problems (see eg. \cite{CCS2,CDS2,HaSc11} 
and the references there). 

Straightforward application of standard
spectral techniques entails, however, the \emph{curse of
dimensionality}: the spectral- or even exponential convergence
rate afforded by analytic parameter dependence is not realized
in computational practice as soon as the number of parameters
is just moderately large. 
High order numerical methods for 
infinite-dimensional problems require, therefore, 
a more refined analysis of analytic parameter dependence
where, for dimension-independent convergence rates, 
the size of the domains of analyticity must increase
with the problem dimension. 

The purpose of the paper is to present recent advances in the analysis of 
higher order \emph{Quasi Monte-Carlo (QMC)} methods,
which were proposed initially in \cite{D08}
(see also \cite{DiPi10}),
from \cite{DKGNS13,DKGSML13}. 
The presented results imply, for a particular
type of analytic parameter dependence encountered for
a large class of operator equations with random coefficients, dimension robust
high order convergence rates, which are only limited by a
certain sparsity measure of the uncertain input.

\section{Affine Parametric Operator Equations}
\label{sec:AffParOpEq}
We present a model setting of affine parametric operator equations,
and their Petrov-Galerkin (PG) discretizations, following the setting
in \cite{DKGSML13}.
We denote by $\bcX$ and $\bcY$ two separable and reflexive Banach spaces
over $\mathbb{R}$ (all results will hold with the obvious modifications
also for spaces over $\mathbb{C}$) with (topological) duals $\bcX'$ and
$\bcY'$, respectively. By $\cL(\bcX,\bcY')$, 
we denote the set of bounded linear operators $A:\bcX \to\bcY'$.
We consider \emph{affine-parametric operator equations}: 
given $f\in \bcY'$, for every $\bsy\in U$ find $u(\bsy)\in \bcX$ 
such that
\begin{equation}\label{eq:main}
  A(\bsy)\, u (\bsy) = f \;.
\end{equation}
%
For such parametrizations, the parametric operator 
$A(\bsy)$ depends on $\bsy$ in an  ``\emph{affine}'' manner:
there exists a sequence $\{ A_j\}_{j\geq 0} \subset \cL(\bcX,\bcY')$ 
such that
\begin{equation}\label{eq:Baffine}
\forall \bsy \in U:\quad 
  A(\bsy) = A_0 + \sum_{j\ge 1} y_j\, A_j \;.
\end{equation}
After possibly rescaling,
we restrict ourselves to the bounded (infinite-dimensional) parameter
domain $U = [-\tfrac{1}{2},\tfrac{1}{2}]^\bbN$.
For every $f\in \bcY'$ and for every $\bsy\in U$, we solve the parametric
operator equation \eqref{eq:main}, where the operator
$A(\bsy)\in\cL(\bcX,\bcY')$ is of affine parameter dependence, see
\eqref{eq:Baffine}. 
We associate with the $A_j$ 
bilinear forms 
$\fa_j(\cdot,\cdot):\bcX\times \bcY \rightarrow \mathbb{R}$ 
via
$$
  \forall v\in \cX,\;w\in \cY:\quad
  \fa_j(v,w) \,=\, {_{\cY'}}\langle  A_j v, w \rangle_{\cY}\;,
  \quad j=0,1,2,\ldots
  \;.
$$
Similarly, for $\bsy\in U$ we associate with the affine-parametric 
operator family $A(\bsy)$ the
parametric bilinear form $\fa(\bsy;\cdot,\cdot): \bcX\times\bcY\to\R$ 
via
\[ 
  \forall v\in \cX,\;w\in \cY:\quad
  \fa(\bsy;v,w) \,=\, {_{\cY'}}\langle A(\bsy) v, w\rangle_{\cY}\;.
\] 
In order for the sum in \eqref{eq:Baffine} to converge, we impose 
\begin{assumption}\label{ass:AssBj}
The sequence $\{ A_j \}_{j\geq 0}\subset \cL(\bcX,\bcY')$ in
\eqref{eq:Baffine} satisfies: 
\begin{enumerate}
\item%
$A_0\in \cL(\bcX,\bcY')$ is boundedly invertible, i.e., there 
exists $\mu_0 > 0$ such that
\begin{equation*} 
 \inf_{0\ne v \in \bcX} \sup_{0\ne w \in \bcY}
 \frac{\fa_0(v,w)}{\| v \|_{\bcX} \|w\|_{\bcY}}
 \ge \mu_0\;,\quad
 \inf_{0\ne w \in \bcY} \sup_{0\ne v \in \bcX}
 \frac{\fa_0(v,w)}{\| v \|_{\bcX} \|w\|_{\bcY}}
 \ge \mu_0
 \;.
\end{equation*}
\item%
The \emph{fluctuation operators} $\{ A_j \}_{j\geq 1}$ 
are small with respect to $A_0$ in the following sense:
there exists a constant $0 < \kappa < 2$ such that
\begin{equation} \label{eq:Bjsmall} 
 \sum_{j\geq 1} \beta_{0,j} \leq \kappa < 2 \;,
 \quad\mbox{where}\quad
 \beta_{0,j} \,:=\, \| A_0^{-1} A_j \|_{\cL(\cX,\cX)}\;,
 \quad j=1,2,\ldots
 \;.
\end{equation}
\end{enumerate}
\end{assumption}
%
%
\begin{theorem}[{cf.~\cite[Theorem 2]{ScMCQMC12}}]
\label{thm:BsigmaInv} Under Assumption~\ref{ass:AssBj}, for every
realization $\bsy\in U$ of the parameter vector, the affine parametric
operator $A(\bsy)$ given by \eqref{eq:Baffine} is boundedly invertible,
uniformly with respect to $\bsy$.
In particular, for every $f \in \bcY'$ and for every $\bsy \in U$, 
the parametric operator equation
\begin{equation} \label{eq:parmOpEq}
 \mbox{find} \quad u(\bsy) \in \bcX:\quad
 \fa(\bsy;u(\bsy), w) \,=\,  {_{\bcY'}}\langle  f,w \rangle_{\bcY}
 \quad
 \forall w \in \bcY
\end{equation}
admits a unique solution $u(\bsy)$ which satisfies
the a-priori estimate
\[ 
 \| u(\bsy) \|_{\bcX}
 \,\leq\,
 \frac{1}{\mu}
 \, \| f \|_{\bcY'}\;, \quad
  \mbox{with}\quad \mu = (1 - \kappa/2)\,\mu_0
\;.
\] 
\end{theorem}
%

\subsection{Single-level and multi-level algorithms}\label{ssec_algorithm}

The Quantity of Interest (QoI) in our study is the expected value of a linear functional 
$G:\cX \rightarrow \bbR$ of the solution $u$,
\begin{equation*}
I(G(u)) = \int_U G(u(\bsy)) \,\mathrm{d} \bsy.
\end{equation*}
In the following we discuss the approximation of the QoI by the algorithm $Q_{N,s}(G(u^h_s))$, where $Q_{N,s}$ is a quadrature rule (QMC rule) and $u^h_s$ is the Petrov-Galerkin (PG) approximation of the dimension truncated problem, which means that the set of parameters $\bsy \in U$ is restricted to $\bsy$ of the form $(y_1, y_2, \ldots, y_s, 0, 0, \ldots)$. The combined error of this \emph{single-level algorithm} can be expressed as
\begin{align}\label{combined_error}
I(G(u)) - &  Q_{N,s}(G(u^h_s)) \nonumber \\ = & \underbrace{I(G(u)) - I(G(u_s)) }_{\mbox{truncation error} } + \underbrace{  I(G(u_s)) - Q_{N,s}(G(u_s)) }_{\mbox{integration error}}  + \underbrace{ Q_{N,s}(G(u_s -u^h_s )) }_{\mbox{PG error}},
\end{align}
where 'PG error' stands for the Petrov-Galerkin discretization error. We discuss the three errors and the necessary background in the subsequent sections.

To reduce the computational cost required to achieve the same error, a
novel \emph{multi-level} algorithm was introduced and analyzed in
\cite{KSS13}. It takes the form
\begin{equation}\label{eq:QL*}
 Q^L_*(G(u)) \,:=\,
 \sum_{\ell=0}^L Q_{s_\ell,N_\ell}(G(u^{h_\ell}_{s_\ell} - u^{h_{\ell-1}}_{s_{\ell-1}}))\;.
\end{equation}
In \cite{KSS13} the authors considered the case where each $Q_{s_\ell,N_\ell}$ is a randomly shifted lattice rule with
$N_\ell$ points in $s_\ell$ dimensions, and where $u^{h_{-1}}_{s_{-1}}:=0$, whereas in \cite{DKGSML13} the authors used an interlaced polynomial lattice rule.

It is well known \cite{DKGNS13} that under some assumptions the Petrov-Galerkin discretization 
error is of the form
\begin{align} \label{eq:Gconvest}
   \left| G(u(\bsy)) - G(u^h(\bsy)) \right|
   &\,\le\, C\, h^{t+t'} \, \| f \|_{\cY'_{t}} \, \| G \|_{\cX'_{t'}}
 \;.
\end{align}

\subsection{Parametric and spatial regularity of solutions}
\label{ssec:anadepsol}
First we establish the regularity of the solution $u(\bsy)$ of the
parametric, variational problem \eqref{eq:parmOpEq} with respect to the
parameter vector $\bsy$. This is important for the analysis of the integration error using a QMC rule satisfying a dimension-independent error bound.

In the following, let $\N_0^\N$ denote the set of
sequences $\bsnu = (\nu_j)_{j\geq 1}$ of non-negative integers $\nu_j$, and
let $|\bsnu| := \sum_{j\geq 1} \nu_j$. For $|\bsnu|<\infty$, we denote the
partial derivative of order $\bsnu$ of $u(\bsy)$ with respect to $\bsy$ by
\[
\partial^\bsnu_\bsy u(\bsy) \,:=\,
\frac{\partial^{|\bsnu|}}{\partial^{\nu_1}_{y_1}\partial^{\nu_2}_{y_2}\cdots}u(\bsy),
\quad \bsy \in U 
\;.
\]
%
\begin{theorem}[cf. \cite{CDS2,KunothCS2011}]
\label{thm:Dsibound}
Under Assumption~\ref{ass:AssBj}, there exists a constant $C_0 > 0$ such
that for every $f\in \bcY'$ and for every $\bsy\in U$, the partial
derivatives of the parametric solution $u(\bsy)$ of the parametric
operator equation \eqref{eq:main} with affine parametric, linear operator
\eqref{eq:Baffine} satisfy the bounds
\begin{equation*} 
\|\partial^\bsnu_\bsy u(\bsy)\|_\bcX
\,\le\,
C_0\, |\bsnu|! \,\bsbeta_0^\bsnu\, \| f\|_{\bcY'}
\quad \mbox{for all } \bsnu \in \N_0^\N \mbox{ with } |\bsnu|<\infty
\;,
\end{equation*}
where $0! :=1$, $\bsbeta_0^\bsnu := \prod_{j\ge 1} \beta_{0,j}^{\nu_j}$,
with $\beta_{0,j}$ as in \eqref{eq:Bjsmall}, and $|\bsnu| = \sum_{j \ge 1}
\nu_j$.
\end{theorem}

\medskip
%
Spatial regularity is in {\em scales of smoothness
spaces} $\{ \cX_t \}_{t \geq 0}$, $\{ \cY_t \}_{t\geq 0}$,
i.e.
\begin{equation*} 
\begin{aligned}
 \cX &= \cX_0 \supset \cX_1 \supset \cX_2 \supset \cdots\;,
 &\cY &= \cY_0 \supset \cY_1 \supset \cY_2 \supset \cdots\;,
 \quad\mbox{and}
 \\
 \cX' &= \cX'_0 \supset \cX'_1 \supset \cX'_2 \supset \cdots\;,
 &\cY' &= \cY'_0 \supset \cY'_1 \supset \cY'_2 \supset \cdots
 \;.
\end{aligned}
\end{equation*}
For self-adjoint operators, usually $\cX_t  = \cY_t$. 
\begin{assumption}[{see \cite[Assumption 2]{DKGSML13}}]\label{ass:XtYt}

There exists $\bar{t}\ge 0$ such that
\begin{enumerate}
\item 
For every $t,t'$ satisfying $0\le t,t'  \le \bar{t}$, we have
\begin{equation}\label{eq:Regul}
 \sup_{\bsy\in U} \| A(\bsy)^{-1} \|_{\cL(\cY'_t, \cX_t)} < \infty
 \quad\mbox{and}\quad
 \sup_{\bsy\in U} \| (A^*(\bsy))^{-1} \|_{\cL(\cX'_{t'}, \cY_{t'})}  < \infty\;.
\end{equation}
%
Moreover, there exist summability exponents $0\le p_0 \le p_t \le
p_{\bar{t}}<1$ such that
\begin{equation} \label{eq:psumpsi0}
 \sum_{j \ge 1} \| A_j \|^{p_t}_{\cL(\cX_t,\cY'_t)} < \infty
\;.
\end{equation}
\item
    Let $\bsu(\bsy) = (A(\bsy))^{-1}f$ and
    $w(\bsy) = (A^*(\bsy))^{-1}G$.
    For $0\le t,t'\le\bar{t}$,
    there exist constants $C_t,C_{t'}>0$
    such that for every  $f\in \cY'_t$ and $G\in\cX'_{t'}$ holds
\begin{equation*} 
  \sup_{\bsy\in U} \|u(\bsy)\|_{\cX_t} \le C_t \|f\|_{\cY'_t}
  \quad\mbox{and}\quad
 \sup_{\bsy\in U} \|w(\bsy)\|_{\cY_{t'}} \le C_{t'} \|G\|_{\cX'_{t'}}\;.
\end{equation*}
Moreover, for every $0\le t\le\bar{t}$ there exists a sequence
$\bsbeta_t = (\beta_{t,j})_{j\geq 1}$ 
satisfying
\begin{equation*} 
 \sum_{j\ge 1} \beta_{t,j}^{p_t} \,<\,\infty\;,
\end{equation*}
such that for
every $0\leq t,t' \leq \bar{t}$ and for every
$\bsnu\in\bbN_0^\bbN$ with $|\bsnu|<\infty$
we have
\begin{align*} 
 \sup_{\bsy \in U} \|\partial^{\bsnu}_{\bsy} u (\bsy) \|_{\cX_t}
 &\,\le\, C_t\, |\bsnu|!\, \bsbeta_t^{\bsnu}\, \|f\|_{\cY_t'},
 \\
 \label{eq:refAd}
 \sup_{\bsy\in U} \| \partial^{\bsnu}_{\bsy} w(\bsy)\|_{\bcY_{t'}}
 &\,\le\, C_{t'}\, |\bsnu|!\, \bsbeta_{t'}^\bsnu\, \| G \|_{\bcX'_{t'}}
\;.
\end{align*}
\item The operators $A_j$ are enumerated so that the sequence
    $\bsbeta_0$ in \eqref{eq:Bjsmall} satisfies
\begin{equation} \label{eq:ordered}
  \beta_{0,1} \ge \beta_{0,2} \ge \cdots \ge \beta_{0,j} \ge \, \cdots\;.
\end{equation}
\end{enumerate}
\end{assumption}
\subsection{Dimension truncation}
\label{sec:dimtrunc}
We truncate the infinite sum in \eqref{eq:Baffine} to $s$ 
terms and
solve the corresponding operator equation \eqref{eq:main}
approximately using Galerkin discretization from
two dense, one-parameter families
$\{\cX^h\}\subset \cX$, $\{\cY^h\}\subset \cY$ of subspaces
of $\cX$ and $\cY$:
for $s\in \bbN$ and $\bsy\in U$, we define
\begin{equation*} 
\fa_{s}(\bsy;v,w) := _{\bcY'}\langle A^{(s)}(\bsy) v, w\rangle_{\bcY},
\quad\text{with}\quad
  A^{(s)}(\bsy) :=  A_0 + \sum_{j=1}^{s} y_j A_j.
\end{equation*}
For $0<h \le h_0$ and $\bsy\in U$,
the dimension truncated PG-solution is defined by
\begin{equation} \label{eq:uhs}
\text{find } u^h_s(\bsy) \in \cX^h: \quad
\fa_s(\bsy;u^h_s(\bsy),w^h) = _{\bcY'}\langle f, w^h \rangle_\bcY \quad
\forall w^h\in \cY^h
\;.
\end{equation}
By choosing $\bsy = (y_1,\ldots,y_s,0,0,\ldots)$,
the PG discretization error bound \eqref{eq:Gconvest} remains 
valid for the dimensionally truncated problem \eqref{eq:uhs}. 

%

\begin{theorem}[{cf.~\cite[Theorem~2.6]{DKGNS13}}] \label{thm:trunc}
Under Assumption~\ref{ass:AssBj},
for every $f\in \cY'$, for every $G\in \cX'$,
for every $\bsy\in U$,
for every $s\in\bbN$ and for every $h>0$,
the variational problem \eqref{eq:uhs}
admits a unique solution $u_s^h(\bsy)$
which satisfies
\begin{equation*} 
  |I(G(u^h))- I(G(u^h_s))|
  \,\le\, C\, 
  \|f\|_{\cY'}\, \|G\|_{\cX'}\,
  \bigg(\sum_{j\ge s+1} \beta_{0,j}\bigg)^2
\end{equation*}
for some constant $C>0$ 
independent of $f$, $G$ and of $s$
where $\beta_{0,j}$ is defined in
\eqref{eq:Bjsmall}. In addition, if 
\eqref{eq:psumpsi0} and \eqref{eq:ordered} hold with $p_0<1$, then
\[
  \sum_{j\ge s+1} \beta_{0,j}
  \,\le\,
  \min\left(\frac{1}{1/p_0-1},1\right)
  \bigg(\sum_{j\ge1} \beta_{0,j}^{p_0} \bigg)^{1/p_0}
  s^{-(1/{p_0}-1)}\;.
\]
\end{theorem}
%


\section{Quasi Monte-Carlo quadrature}
\label{sec:QMC}

In \cite{KSS12}, Quasi-Monte Carlo rules of the form $Q_{N,s}(G(u^h_s)) = \frac{1}{N} \sum_{n=0}^{N-1} G(u^h_s(\bsy_n - \mathbf{ \tfrac{1}{2} } ))$, where $\bsy_n \in [0,1]^s$, have been used to approximate the dimension truncated integral $I(G(u^h_s))$ (see also \cite{KSS11}). The rules considered therein are so-called randomly shifted lattice rules. Using so-called ``\emph{product and order-dependent (POD) weights}'' a convergence rate of order $\calO(N^{-\min(1/p_0-1, 1-\delta)})$, for any $\delta > 0$, was shown.

Noting that the integrand is actually analytic, the authors of \cite{DKGNS13} used \emph{interlaced polynomial lattice rules}, 
as introduced in \cite{DiGo12} (which are a special type of higher order digital net \cite{D08}), 
to obtain improved rates of convergence. 
The rules can be constructed using the fast component-by-component approach of \cite{NC06a}. 
A new function space setting was introduced in \cite{DKGNS13} which uses Banach 
spaces and \emph{smoothness driven product and order dependent (SPOD) weights}.

\begin{theorem} [{cf. \cite[Theorem~3.1]{DKGNS13}}] \label{thm:main}
Let $s\ge 1$ and $N = b^m$ for $m\ge 1$ and prime $b$. Let $\bsgamma =
(\gamma_j)_{j\ge 1}$ be a sequence of positive numbers, let $\overline{\bsgamma}_s =
(\gamma_j)_{1\le j \le s}$, and assume that
\begin{equation*} 
  \exists\, 0<p\le 1 : \quad \sum_{j=1}^\infty \gamma_j^p < \infty\;.
\end{equation*}
Define
Suppose we have an integrand $F(\bsy)$ whose partial derivatives satisfy
\begin{equation*} 
 \forall\, \bsnu \in \{0, 1, \ldots, \alpha\}^s: \quad
 | (\partial^{\bsnu}_\bsy F)(\bsy)| \,\le\, c\, |\bsnu|!\, \bar{\bsgamma}_s^{\bsnu}
\end{equation*}
for some constant $c>0$.
Then, an interlaced polynomial lattice rule of
order $\alpha$ with $N$ points can be constructed using a fast
component-by-component algorithm, with cost
$\calO(\alpha\,s\, N\log N + \alpha^2\,s^2 N)$ operations,
such that
\[
  |I_s(F) - Q_{N,s}(F)|
  \,\le\, C_{\alpha,\bsgamma,b,p}\, N^{-1/p} \;,
\]
where $C_{\alpha,\bsgamma,b,p} < \infty$ is a constant independent of $s$
and $N$.
\end{theorem}


\section{Combined error bound}

In the case of the single level algorithm, the combined error \eqref{combined_error} satisfies the following theorem.
\begin{theorem}[{cf. \cite[Theorem~4.1]{DKGNS13}}]
Under Assumption~\ref{ass:AssBj} and conditions \eqref{eq:Regul},
$G \in \cX'_{t'}$ and \eqref{eq:ordered}, the integration error
using an interlaced polynomial lattice rule of order $\alpha = \lfloor 1/p_0 \rfloor
+ 1$ with $N=b^m$ points (with $b$ prime) in $s$ dimensions, combined with
a Petrov-Galerkin method in the domain $D$ with one common subspace
$\cX^h$ with $M_h=\operatorname{dim}(\cX^h)$ degrees of freedom and 
with linear cost $\calO(M_h)$, satisfies
\[
  |I(G(u)) - Q_{N,s}(G(u^h_s)) | \,\le\, \calO\left(s^{-2(1/p_0-1)} + N^{-1/p_0} + h^{t+t'}\right)
   \;,
\]
where the constant is independent of $s$, $h$ and $N$.
\end{theorem}
The multi-level algorithm additionally requires the Assumptions~\ref{ass:XtYt}. 
The corresponding combined error bound using interlaced polynomial lattice rules 
is of the form (see \cite[Theorem~3.4]{DKGSML13})
\[
 |I(G(u)) - Q^L_{*}(G(u^h_s)) | 
\,\le\, 
\calO\left(s_L^{-2(1/p_0-1)} + h_L^{t+t'} + \sum_{\ell=0}^L N_\ell^{-1/p_t} 
\left(s_{\ell-1}^{-(1/p_0 - 1/p_t)} + h_{\ell-1}^{t+t'} \right) \right)
\;.
\]
The parameters $s_\ell$ and $N_\ell$ in \eqref{eq:QL*} can be optimized using a 
Lagrange multiplier argument \cite{KSS13,DKGSML13}, which, in most cases, 
yields an improvement compared to the single-level algorithm.
\begin{acknowledgement}
J. D. is the recipient of an Australian Research Council Queen Elizabeth II Fellowship 
(project number DP1097023). QLG was supported partially by the ARC Discovery Grant DP120101816.
The work of CS was supported in part by 
European Research Council AdG grant STAHDPDE 247277, and the
Swiss National Science Foundation.
\end{acknowledgement}

\end{document}